\documentstyle[12pt]{article}

\catcode `\@=11
\@addtoreset{equation}{section}
\def\theequation{\arabic{section}.\arabic{equation}}
\catcode `\@=12

\newcommand{\be}{\begin{equation}}
\newcommand{\en}{\end{equation}}
\newcommand{\bea}{\begin{eqnarray}}
\newcommand{\ena}{\end{eqnarray}}
\newcommand{\beano}{\begin{eqnarray*}}
\newcommand{\enano}{\end{eqnarray*}}
\newcommand{\bee}{\begin{enumerate}}
\newcommand{\ene}{\end{enumerate}}

\newcommand{\R}{R \!\!\!\! R}
\newcommand{\N}{N \!\!\!\!\! N}
\newcommand{\Z}{Z \!\!\!\!\! Z}
\newcommand{\Hil}{{\cal H}}

\newcommand{\Id}{1\!\!1}
\newcommand{\F}{{\cal F}}
\newcommand{\Lc}{{\cal L}}
\newcommand{\D}{{\cal D}}
\newcommand{\C}{{\cal C}}
\newcommand{\E}{{\cal E}}
\newcommand{\A}{{\cal A}}

\begin{document}

\thispagestyle{empty}
 
\vspace*{1cm}

\begin{center}
{\Large \bf Applications of Topological *-Algebras of Unbounded
Operators}   \vspace{2cm}\\

{\large F. Bagarello}
%\footnote[1]{ Dipartimento di Matematica ed Applicazioni, 
%Fac. Ingegneria, Universit\`a di Palermo, I-90128  Palermo, Italy}  
\vspace{3mm}\\
  Dipartimento di Matematica ed Applicazioni, 
Fac. Ingegneria, Universit\`a di Palermo, I-90128  Palermo, Italy\\
e-mail: Bagarello@ipamat.math.unipa.it
\vspace{2mm}\\
\end{center}

\vspace*{2cm}

\begin{abstract}
\noindent 
In this paper we discuss some physical applications of topological
*-algebras of unbounded operators. Our first example is
 a simple system of free bosons. Then we analyze different models
which are related to this one. We also discuss the time evolution of two
interacting models of matter  and bosons. We show that for all these
systems it is possible to build up a common framework where the
thermodynamical limit of the algebraic dynamics can be conveniently
studied and obtained.   \end{abstract}

\vspace{2cm}


\vfill

\newpage

% Section 1
\section{Introduction}

Since when for the first time the algebraic approach to the
Quantum Field Theory, \cite{hk}, has made its appearance, many new
results have been found and discussed in the literature showing that the
structure of $C^*$-algebras is not rich enough in many relevant
physical cases. Difficulties arise already in ordinary quantum
mechanics, since the commutation rule $[x,p]=i$ implies that not both the
operators $x$ and $p$ can be bounded as operators on $\Lc^2(\R)$. Much
more involved is the situation for models with infinite degrees of
freedom. These difficulties have given rise to a wide literature related
to the subject of unbounded operators: the Lassner's quasi *-algebras,
\cite{las}, the partial *-algebras introduced by Antoine and Karwowsky,
\cite{ant}, the $CQ^*$-algebras recently introduced by Trapani and the
author, \cite{btcq}, etc.. . Many of these results are collected and
analyzed in particular in two books: the first, essentially
mathematically oriented, was written in 1989 by Schm\" {u}dgen,
\cite{sch}; the second, more (but not to much!) thought for a mathematical physicist
reader, was conceived one year later by Dubin and Hennings, \cite{dub}.

However, it is impossible not to notice that all this work, with few
exceptions, is essentially concerned with purely mathematical aspects
of the problem, while physical applications are missing, for the time
being. In other words, in the course of the years many powerful tools
have been proposed, constructed and refined, but not many real physical
problems have been solved or even settled using these tools. The present
paper try to make a step to fill this gap. We use topological *-algebras
of unbounded operators, of the kind discussed first by Lassner,
\cite{las}, to discuss the existence of the algebraic dynamics of
different physical models. We begin with an apparently innocuous model of
free bosons, which already presents the difficulty of dealing with
unbounded operators, and we show how this model can be regularized, using
a certain cutoff, and then we prove the existence of the algebraic
dynamics of the model when the cutoff is {\em cleverly } removed. This is
the content of the first part of Section 3. In the second part of the
same Section we analyze different models which can be related again to
the free bosons using simple transformations. In Section 4 we consider
first an interacting model of single mode bosons interacting with the
matter, mimicked by spin observables localized at different lattice sites,
and we show again the possibility of regularizing the hamiltonian and to
find a limit of the time evolution in a natural topology, again related
to the topologies first introduced by Lassner, \cite{las}. Then
we extend our procedure to a generalized version of the same
model, where different kind of bosons are considered. Section 5 contains
the outcome of the paper. We begin, in the following Section, fixing the
mathematical structure of our approach. Some of the mathematical details
will be discussed also in the Appendices.

Before starting a remark is in order: the topological structure of the
algebra of the unbounded operators is not always essential. In fact,
depending on the approach used for removing the cutoff, it can be
sufficient to work in an algebra without topology. This has been done,
for instance, in \cite{alsew, bagsew}, considering weak limits of the
time evolution of local and delocalized observables.

% Section 2
\section{The Mathematical Framework}
In this paper we will 
consider mainly a single mode boson oscillator. Here $a$ and $a^\dagger$
are respectively the annihilation and the creation operators, satisfying 
the CCR \be
[a,a^\dagger]=\Id.
\label{21} 
\en

We  
introduce in a canonical way the Hilbert space of the radiation: we
consider a vector $\Phi_0$ which is annihiled by $a$, and then we build up
the linear span of the vectors $\{(a^\dagger)^n\Phi_0, n\in \N\}$,
$\D_o$. Taking the completion of this set, we get the Hilbert space of
the bosons, $\Hil$, which is, of course, a Fock space.

Let now $N_o\equiv a^\dagger a$. This is a symmetric operator on $\D_o$
which can be uniquely extended to a self-adjoint operator $N$. Let then
$D(N^k)$ be the domain of the operator $N^k$, $k\in \N$. We also
introduce a set $\D$ as the domain of all the powers of $N$:  \be
\D \equiv D^\infty(N) = \cap_{k\geq 0} D(N^k).
\label{22}
\en
This set is dense in $\Hil$. Starting from $\D$ we can define, following
Lassner, \cite{las}, the *-algebra $\Lc^+(\D)$ as the set of all
the closable operators defined on $\D$ which map, together with their
adjoints, $\D$ into itself. It is clear that all the powers of $a$ and
$a^\dagger$ belong to such a set.

In \cite{las} the topological structures of both  $\D$ and $\Lc^+(\D)$
are discussed in details, and the role played by the number operator $N$
in defining all these seminorms is discussed. In particular, in $\D$ the 
seminorms are 
\be \phi \in \D \rightarrow \|\phi\|_n\equiv \|N^n\phi\|,
\label{23} \en
where $n$ is a natural integer while $\|\, \|$ is the norm in $\Hil$. The
topology in $\Lc^+(\D)$ is introduced in the following way. We start
defining the set  $\C$ of all the positive, bounded and continuous
functions $f(x)$ on $\R_+$, which are decreasing faster than any inverse
power. The seminorms on  $\Lc^+(\D)$ are labeled by functions in $\C$ and
by the integers $\N$. We have 
\be
X\in \Lc^+(\D) \rightarrow \|X\|^{f,k} \equiv
\max\left\{\|f(N)XN^k\|,\|N^kXf(N)\|\right\}. 
\label{24}
\en
Here $\|\,\|$ is the usual norm in $B(\Hil)$. We use for this norm the
same notation as in equation (\ref{23}) since no confusion can arise.
Incidentally, we see that a possible way to look at
 $\Lc^+(\D)$ is as the set of operators such that both $f(N)XN^k$ and
$N^kXf(N)$ are bounded for all integers $k$ and for all functions $f$ in
$\C$. Moreover, we observe that an easy consequence  of the definition
(\ref{24}) is the following invariance property:
$\|X\|^{f,k}=\|X^\dagger\|^{f,k}$. We call $\tau$ the topology on
$\Lc^+(\D)$ defined by the seminorms in (\ref{24}). In \cite{las} it has
been proven that  $\Lc^+(\D)[\tau]$ is a locally convex complete
topological *-algebra.

In the following Sections a useful role will be played by some
subspaces of $\Hil$ which we now introduce. We start defining the
one-dimensional spaces $\E_l$ as the set of all the vectors which are
proportional to $(a^\dagger)^l\Phi_0$. Then we introduce the
$(L+1)$-dimensional space $\F_L\equiv
\E_0\oplus\E_1\oplus.....\oplus\E_L$, as the direct sum of the first
$L+1$ spaces $\E_l$.

These spaces are in one-to-one correspondence with some projection
operators. Let $N=\sum_{l=0}^{\infty}l\Pi_l$ be the spectral decomposition
of the number operator $N$. The operators $\Pi_l$ are 
projection operators, as well as the operators 
$Q_L=\sum_{l=0}^{L}\Pi_l$. The following properties are, therefore,
consequences of their nature: 
\be \Pi_k \Pi_l=\delta_{kl}\Pi_l, \hspace{1cm}
\Pi_k^\dagger=\Pi_k,  \label{25}
\en
and
\be
Q_L Q_M=Q_L, \hspace{.5cm} \mbox{if } L\leq M, \hspace{1cm}
Q_L^\dagger=Q_L. 
\label{26}
\en
Obviously $\Pi_k:\Hil \rightarrow \E_k$, while $Q_L:\Hil \rightarrow
\F_L$. Equations ((\ref{25}),(\ref{26})) show also that the spaces $\E_k$
are mutually orthogonal, while the $\F_L$ are not, properties which also
can be derived from the original definitions of the spaces. Of course, the
following inclusions of spaces hold:
\be
\E_l\subset \F_L \subset \F_{L+1} \subset \D \subset \Hil,
\label{26bis}
\en
whenever $l\leq L$.

In this paper we will not make use of another possibility which sometimes
can be of a certain interest, that is the one considering, instead of
$\Lc^+(\D)$, the set $\Lc(\D, \D')$ of the linear continuous maps from
$\D$ into its dual $\D'$. Also in this set a topology can be defined, see
\cite{las}, which makes it complete, but what makes this space less
interesting than $\Lc^+(\D)$ for our purposes is the lack of a convenient
algebraic structure. In other words, while it is well posed the problem
of finding the time evolution for the product $AB$ when both $A$ and $B$
belong to $\Lc^+(\D)$, this problem cannot be discussed in general if $A$
and $B$ are both elements of $\Lc(\D, \D')$. As it is well known, in this
case only partial multiplications could be defined, within elements in 
$\Lc^+(\D)$ and elements of $\Lc(\D, \D')$. A rather complete review on
this and related subjects was written by C. Trapani, \cite{tra}.

\vspace{3mm}

Some final comments concerning the seminorms in
(\ref{24}). The first remark is that the two contributions in the
definition are exactly of the same form. The estimate of $\|f(N)XN^k\|$ is
very similar to the estimate of $\|N^kXf(N)\|$. This is why in the
following we will only consider the first of these contributions,
identifying in this way $\|X\|^{f,k}$ simply with 
$\|f(N)XN^k\|$. 

With this in mind, we came to the second remark, which is a consequence
of the spectral decomposition for $N$. Using $N=\sum_{l=0}^\infty l
\Pi_l$, we see that the seminorms can be written as follows: \be
X\in \Lc^+(\D) \longrightarrow \|X\|^{f,k} = \sum_{l,s=0}^{\infty}
f(l)s^k\|\Pi_lX\Pi_s\|.
\label{27}
\en
From now on we will denote by {\em physical } the topology $\tau$ generated
by these seminorms, following with some freedom the notation usually
adopted in the literature. 

\vspace{5mm}

Since in Section 4 we will discuss a model which, besides a boson
component contains also spin variables, we introduce here some more
relevant spaces. We refer to \cite{bm} for further details.

Let $\Z^d$ be a d-dimensional infinite lattice, $V\subset \Z^d$, and
$|V|$ the number of the points of $V$. We call $\A_V$ the $C^*$-algebra
generated by the spin operators $\sigma_\alpha^i$, $i\in V$ and
$\alpha=x,y,z$, and $\A_o$ the norm closure of $\cup_V \A_V$.

We call {\em relevant } a state $\omega$ over $\A_o$ if, denoting by
$\Hil_\omega$ the Hilbert space defined by the GNS construction on $\A_o$
and $\omega$, and by $\Psi_\omega$ the vector which represents $\omega$
in $\Hil_\omega$ ($\omega(A)=<\Psi_\omega, \pi_\omega(A)\Psi_\omega>$, for
all $A\in \A_o$), then $\Psi_\omega$ belongs to the following set $\F$:
\be
\F =\left\{ \Psi \in \Hil_\omega \: :
\lim_{|V|,\infty} \frac{1}{|V|}\sum_{p\in V} {\sigma_3^p}\Psi =
\sigma_3^\infty \Psi,  \:\: \|\sigma_3^\infty\|\leq1 \right\},
\label{28}
\en
where $\sigma_3^\infty$ belongs to the center of the algebra $\A_o$ and
$\pi_\omega$ is the canonical representation for the spin algebra,
\cite{bm}.

The relevance of this set has been discussed in reference \cite{bm},
where it has been proven, among other things, that all the powers of
$\sigma_3^V$ converge in the $\F$-strong topology (the strong
topology 'restricted' to those vectors which belong to $\F$), as well as
the analytic functions of $\sigma_3^V$. This reflects the well known fact
that $\sigma_3^V$ is not norm converging for $|V|$ increasing, \cite{thi}.
The topology is therefore defined by the following seminorms:
\be
X\in \A_o \simeq B(\Hil_{spin}) \rightarrow \|X\|^\Psi \equiv \|X\Psi\|,
\label{29}
\en
where $\Psi \in \F$, $X$ is identified with its 'canonical'
representative in $B(\Hil_{spin})$, $\Hil_{spin}$ being the infinite
tensor product of two-dimensional complex spaces ${\bf C}^2_i$, $i\in
\Z^d$. In this way $\A_o$ is canonically identified with its
representative $B(\Hil_{spin})$. This identification will be used all
throughout this paper in order to simplify the notation.

% Section 3
\section{The Free Bosons and Related Models}

%3.1 
\subsection{\hspace{-6mm}. The free model}

In this Subsection we propose a possible rigorous
algebraic approach to the problem of the time evolution of a
single mode of radiation described by the following free hamiltonian:
\be
H=a^\dagger a,
\label{31}
\en
where $a$ and $a^\dagger$ are the boson operators introduced in 
 Section 2. In all this paper we will follow the canonical strategy 
present everywhere in the literature, whenever a certain rigor is
required: first we regularize the hamiltonian, by means of a certain
cutoff, then the equation of motion are obtained and solved keeping fixed
this cutoff which is finally removed. If all these steps can be
performed, we define the dynamics for the 'infinite' model as the limit
of these cutoffed dynamics.

The first non trivial problem is the way in which this cutoff should be
implemented. As it is obvious from the expression of $H$, which is nothing
but the number operator $N$, $H$ is unbounded. Nevertheless,
if we consider $N$ acting not on the whole $\Hil$ but only on some
(cleverly chosen) subspaces of the Hilbert space, the unboundedness can be
controlled. In fact, it is easily seen that $\|N\varphi\|<c_\varphi$
whenever $\varphi$ is taken in any of the spaces $\E_l$ or $\F_L$, where
$c_\varphi$ is a positive constant depending on $\varphi$ (and therefore
on $l$ and $L$). This fact suggests to consider the following
regularization:  \be H_L=Q_L H Q_L
\label{32}
\en
which preserve the hermitianness of the hamiltonian and satisfies the
following property: $\forall \varphi \in \Hil \Rightarrow
H_L\varphi=Q_LH\varphi_L$, where $\varphi_L \in \F_L$. Therefore, $\forall
\varphi \in \Hil,\|H_L\varphi\|<c_{\varphi_L}$. However, this
regularization has an unpleasant drawback which is evident already for
the model we will discuss in some details in Section 4, and which is
described by the following  hamiltonian,   $$
H_V=\frac{1}{|V|}\sum_{i,j\in V}\sigma_3^i\sigma_3^j+a^\dagger a +\gamma
(a+a^\dagger)\sigma_3^V,  $$
 where $\sigma_3^V=\frac{1}{|V|}\sum_{i\in V}\sigma_3^i$. In fact, it is
clear that, after a regularization procedure as in (\ref{32}), the
following term would appear: $\frac{1}{|V|}\sum_{i,j\in
V}\sigma_3^i\sigma_3^jQ_L$. This contribution is not certainly natural,
and it is also technically difficult to be properly considered, since the
bosonic nature of the problem is now mixed to its spin component, already
in the free spin hamiltonian.

All these considerations suggest the use of another regularization,
defined starting from the boson operators themselves. We perform the
following substitutions: 
\be
a \rightarrow a_L\equiv Q_LaQ_L, \hspace{1cm} a^\dagger \rightarrow
a^\dagger_L\equiv Q_La^\dagger Q_L,
\label{33}
\en
and then we define, for the model in (\ref{31}),
\be
H_L = a^\dagger_L a_L.
\label{34}
\en
This procedure can be naturally generalized: if $H$ is a given function
$\chi$ of certain non bosonic (bounded) operators $X_1,..,X_n$, as well as
of $a$ and $a^\dagger$, eventually already depending on some cutoff as in
the hamiltonian $H_V$ above, we can consider the {\em regularized}
hamiltonian $H_L$ as the same function of the same variables $X_1,..,X_n$
and of the regularized boson operators $a_L$ and of $a^\dagger_L$:
\be
H=\chi(X_1,..,X_n,a,a^\dagger) \rightarrow
H_L=\chi(X_1,..,X_n,a_L,a_L^\dagger).
\label{34bis}
\en

\vspace{4mm}

{\bf Remark.}-- Using the results discussed in Appendix A it is easy to
show that for the free hamiltonian $H=N$ the above regularizations
coincide. In fact, by equations (\ref{a2}) and (\ref{26}), we have
$$
a^\dagger_La_L=Q_La^\dagger Q_LQ_LaQ_L=Q_La^\dagger Q_LaQ_L=Q_La^\dagger
aQ_{L+1}Q_L=Q_LNQ_L.
$$
This is not true for general models, as the spin-bosons hamiltonian above
explicitly shows.

\vspace{4mm}

Now that we have chosen the regularization, we come to the main problem 
of this Section, which is to prove the existence of the limit, for $L$
going to infinity, of the following algebraic dynamics: $$
\alpha_L^t(X) = e^{iH_Lt}Xe^{-iH_Lt},
$$
where $X$ is an element of $\Lc^+(\D)$. The topology in which this limit
will be considered is the physical one defined in the previous
Section, which makes $\Lc^+(\D)$ a complete topological *-algebra. We state
the result as a Proposition.

\vspace{4mm}

 {\bf Proposition 3.1.} -- The limits of $\alpha_L^t(a^n)$ and
$\alpha_L^t((a^\dagger)^n)$ exist in $\Lc^+(\D)[\tau]$   for
all natural $n$.

\vspace{3mm}

{\bf Proof}

We prove the Proposition only for the anniihilation operator $a$. The proof
for $a^\dagger$ is essentially the same. The technique we use is an
induction on $n$. We start therefore proving the statement for $n=1$. 

The first step is to observe that:
\be
[H_L,a]_m=(L\Pi_L-Q_{L-1})^ma=(L^m\Pi_L+(-1)^mQ_{L-1})a.
\label{35}
\en
Again, the proof of these formulas goes on using induction on $m$,
and is a consequence of the commutation relations proved in Appendix A.
The equality  $(L\Pi_L-Q_{L-1})^m=(L^m\Pi_L+(-1)^mQ_{L-1})$ is a
consequence of the properties of the projection operators $\Pi_l$ and
$Q_L$. In particular, it is crucial that $\Pi_L
Q_{L-1}=0$, due to equation (\ref{25}) and to the definition of $Q_L$.

By means of (\ref{35}) we can prove that
\be
\alpha_L^t(a)=e^{iH_Lt}ae^{-iH_Lt}=\left(1-Q_L+e^{itL}\Pi_L+e^{-it}
Q_{L-1}\right) a.
\label{36}
\en
In fact, since $H_L$ is bounded, we have
$$
\alpha_L^t(a)= \sum_{m=0}^\infty\frac{(it)^m}{m!}[H_L,a]_m=
\sum_{m=0}^\infty\frac{(it)^m(L^m\Pi_L+(-1)^mQ_{L-1})a
}{m!},
$$
so that equation (\ref{36}) can be obtained simply by resumming this
series.

Now we are ready to prove the $\tau$-Cauchy nature of $\alpha_L^t(a)$.
Using the seminorms in (\ref{27}) we have
$$
Q_{L,M}(t)\equiv \|\alpha_L^t(a)-\alpha_M^t(a)\|^{f,k}=
\sum_{l,s=0}^{\infty}f(l) s^k \|\Pi_l(\alpha_L^t(a)-\alpha_M^t(a))\Pi_s\|.
$$ 
We are not expliciting the dependence of $Q_{L,M}(t)$ from the
seminorm, that is from $f$ and $k$. Using the explicit expression
(\ref{36}) for $\alpha_L^t(a)$ we get the following inequality: \beano
Q_{L,M}(t) &\leq &\sum_{l,s=0}^{\infty}f(l) s^k
{\LARGE [}
\|\Pi_l(Q_M-Q_L)a\Pi_s\|+|e^{-it}|\|\Pi_l(Q_{L-1}-Q_{M-1})a\Pi_s\|+\\
&+&|e^{itL}|\|\Pi_l\Pi_La\Pi_s\|+|e^{itM}|\|\Pi_l\Pi_Ma\Pi_s\|{\LARGE
]}. \enano
To conclude we use equation (\ref{25}) and formula (\ref{a4}):
\be
\|\Pi_la\Pi_s\|=\delta_{l,s-1}\sqrt{s},
\label{37}
\en
so that we get the following estimate for
$Q_{L,M}(t)$: \be
Q_{L,M}(t) \leq 2\sum_{s=M+1}^{L}f(s-1)s^{k+1/2} \rightarrow 0,
\label{38}
\en
when $L$ and $M$ both diverge, due to the nature of the functions in the
set $\C$. Since $\Lc^+(\D)$ is $\tau$-complete, the above result implies
the existence of the limit of $\alpha_L^t(a)$ in $\Lc^+(\D)[\tau]$.

\vspace{4mm}

{\bf Remark:}-- As we have already discussed in Section 2, we have
focused the attention on the seminorm $\|f(N) ..N^k\|$. It is clear from
the above computations, how to extend the same estimates to
$\|N^k..f(N)\|$, which goes to zero essentially for the same reasons. This
implies, due the invariance of the 'real' physical topology, given by the
seminorms in (\ref{24}), with respect to the adjoint, that also
$\alpha_L^t(a^\dagger)$ is $\tau$-Cauchy and therefore $\tau$-convergent
in  $\Lc^+(\D)[\tau]$.

\vspace{4mm}

At this point we can prove the second step of the induction. Hence, we
suppose that $\alpha_L^t(a^{n-1})$ is $\tau$-Cauchy and we use this
hypothesis to prove that also $\alpha_L^t(a^{n})$ is $\tau$-Cauchy. The
proof goes like follows:
\bea
Q_{L,M}^n(t)&\equiv&\|\alpha_L^t(a^n)-\alpha_M^t(a^n)\|^{f,k}\leq
\|\alpha_L^t(a)((\alpha_L^t(a))^{n-1}-(\alpha_M^t(a))^{n-1})\|^{f,k}+
\nonumber\\
&+&\|(\alpha_L^t(a)-\alpha_M^t(a))(\alpha_M^t(a))^{n-1}\|^{f,k}
\label{38bis}
\ena
The first contribution can be estimated by the seminorm
$\|(\alpha_L^t(a^{n-1})-\alpha_M^t(a^{n-1}))\|^{xf,k}$, which we know by
the induction hypothesis to be converging to zero since, if $f(x)\in \C$,
then $xf(x)\in \C$ as well. Without going too much into details, let us
first observe that, calling $F_L(t)\equiv 1-Q_L+e^{itL}\Pi_L+e^{-it}
Q_{L-1}$, then $\alpha_L^t(a)=F_L(t)a$, that $\|F_L(t)\|\leq 3$ and that
$F_L(t)$ commutes with $\Pi_l$ (and $Q_L$). Therefore the first
contribution above can be estimated with
$$
3\sum_{l,s=0}^\infty f(l)s^k \|a\Pi_{l+1}\|\|\Pi_{L+1}(\alpha_L^t(a))^{n-1}
-(\alpha_M^t(a))^{n-1})\Pi_s\|,
$$
so that it is easy to verify our claim above resumming over $s$ and
$l$.

The second contribution in (\ref{38bis}) can be estimated with similar
techniques: again we use formula (\ref{27}) for the seminorms and, using
the commutativity of $F_L(t)$ and $\Pi_l$ together with formula
(\ref{a1}), we get
\beano
&& \|\Pi_l(\alpha_L^t(a)-\alpha_M^t(a))(F_M(t)a)^{n-1}\Pi_s\|= \\
&&=\|\Pi_l(\alpha_L^t(a)-\alpha_M^t(a))\Pi_{s-(n-1)}F_M(t)a\Pi_{s-(n-2)}F_M(t)
a.....\Pi_{s-1} F_M(t)a\Pi_s\|\leq \\
&&\leq
\|\Pi_l(\alpha_L^t(a)-\alpha_M^t(a))\Pi_{s-n+1}\|\: 3^{n-1}\sqrt{(s-n+2) 
(s-n+3)\cdot \cdot \cdot (s-2)(s-1)s} \leq \\
&&\leq (3\sqrt{s})^{n-1}
\|\Pi_l(\alpha_L^t(a)-\alpha_M^t(a))\Pi_{s-n+1}\|.
\enano
Here we have used also the estimate $\|a\Pi_l\|\leq \sqrt{l}$. After some
more easy steps, we finally find $$
\|(\alpha_L^t(a)-\alpha_M^t(a))(\alpha_M^t(a))^{n-1}\|^{f,k}\leq
2\sum_{l=M}^Lf(l)(l+n)^k\sqrt{l+1}\, (3\sqrt{l+n})^{n-1}, $$
which, again, goes to zero for any fixed $n$ when $L$ and $M$ go to
infinity, due to the decay properties of the function $f\in \C$.

\hfill $\Box$

%3.2 
\subsection{\hspace{-6mm}. Variations on the Same Theme}

In this subsection we discuss briefly two more models which can be easily
translated into the same model discussed previously.

The first model is simply a linear perturbation of the original one, and
is described by the following hamiltonian:
$$
H=a^\dagger a+\gamma (a+a^\dagger),
$$
where $\gamma$ is real. As it is well known, the above model
can be rewritten as a free model in terms of the different boson operators
$$
b\equiv a+\gamma.
$$
The operator $b$ and its adjoint $b^\dagger$ satisfy the canonical
commutation relation, as well as $a$ and $a^\dagger$. Furthermore, in
terms of $b$, we can write $H=b^\dagger b-\gamma^2$. This suggests how to
proceed: we introduce the number operator for $b$, $N_b=b^\dagger b$, and
then consider its spectral decomposition, instead of the one for
$N_a=a^\dagger a$, and proceed as in the previous subsection. Every step
can be repeated also in this case, so that even now it is possible to
define rigorously the model first introducing, and then removing, the
cutoff $b\rightarrow b_L=Q_LbQ_L$.

\vspace{4mm}

The second model we are going to analyze presents some more interesting
aspects. It is a model of two modes of radiation coupled to each other. We
assume for it the following hamiltonian  $$
H=a^\dagger a + b^\dagger b +(a^\dagger b+b^\dagger a).
$$
Here $a$ and $b$ both satisfy the CCR and are completely decoupled:
$[a^\sharp,b^\sharp]=0$, where with $X^\sharp$ we denote $X$ or
$X^\dagger$.

Again, the hamiltonian can be mapped into something very similar to a
number operator. We introduce the linear combinations \be
A=\frac{a+b}{\sqrt{2}}, \quad \quad B=\frac{a-b}{\sqrt{2}}.
\label{39}
\en
With these definitions we get $H=2A^\dagger A$. The commutation
relations are  
$$
[A,A^\dagger]=1, \quad [B,B^\dagger]=1, \quad [A^\sharp, B^\sharp]=0.
$$
We can consider now the spectral decomposition for the operator
$N_A=A^\dagger A=\sum_{l=0}^\infty l\Pi_l$, and from the above
commutation relations it is clear that $\Pi_l$ also commutes with $B$ and
$B^\dagger$. At this point, with the usual regularization $H
\rightarrow H_L=2A_L^\dagger A_L$, where $A_L=Q_LAQ_L$ and
$Q_L=\Pi_0+\Pi_1+....+\Pi_L$, the limit for $L\rightarrow \infty$ of
$\alpha_L^t(A^n)$ exists in the topology $\tau$ defined as it is
described in Section 2 starting, this time, with $N_A$. Furthermore, it is
easily checked that $\alpha_L^t(B^n)=B^n$ for all integers $n$, so that
its limit for $L$ diverging exists as well. Therefore, inverting equations
(\ref{39}), we find that also $\alpha_L^t(a^n)$ and $\alpha_L^t(b^n)$ do
converge in $\Lc(\D)[\tau]$, for all integers $n$.

% Section 4
\section{Two Interacting Models}

%4.1 
\subsection{\hspace{-6mm}. Single mode bosons}

In this subsection we will show how to use the framework introduced in
Section 2 to discuss the existence of the algebraic dynamics for the
model described by the finite-volume hamiltonian already introduced in
Section 3: 
\be
H'_V=\frac{J}{|V|}\sum_{i,j\in V}\sigma_3^i\sigma_3^j+a^\dagger a +\gamma
(a+a^\dagger)\sigma_3^V,
\label{41}
\en
where $\sigma_3^V=\frac{1}{|V|}\sum_{i\in V}\sigma_3^i$. This is a
prototype model of an interaction between the matter, considered
as a family of two levels atoms, and conveniently mimicked using
spin matrices, with a single boson mode. 

We take $\gamma$ as a 'small' real parameter, in a sense which will
appear clear in the following, and $V$ as a subset of the infinite lattice
$\Z^d$. We remind that the Pauli spin matrices satisfy the following
commutation rules: \be
[\sigma_\alpha^i,\sigma_\beta^j]=2i\epsilon_{\alpha \beta \gamma}
\delta_{i,j} \sigma_\gamma^i. \label{41bis}
\en
A consequence of this commutation relation, together with
the definition of $\sigma_3^V$, is
\be
[\sigma_3^V,\sigma_\beta^j]=\frac{2i\epsilon_{3 \beta \gamma}}{|V|}
\sigma_\gamma^i,
\label{41ter}
\en
which shows already that, in the limit $|V|\rightarrow \infty$,
$\sigma_3^V$ commutes with all the other observables.

It is worthwhile to focus reader's attention to the volume cutoff in
$H_V'$. The necessity of this cutoff in a rigorous treatment of the
mean field version of the BCS model for the high temperature
superconductivity has been discussed, for instance, in \cite{thi}.

Assuming that the perturbation $\gamma (a+a^\dagger)\sigma_3^V$ is small
compared with the energy difference between the energy levels of the
harmonic oscillator, we can consider 'frozen' the boson kinetic
contribution in $H'_V$. Therefore we take
 \be 
H_V=\frac{J}{|V|}\sum_{i,j\in
V}\sigma_3^i\sigma_3^j+\gamma (a+a^\dagger)\sigma_3^V.
\label{42}
\en
What we have done is the same kind of approximation which
leads to the hamiltonian for the Fractional Quantum Hall Effect starting
with an interacting two-dimensional electron gas in presence of a strong
magnetic field. In that system, we consider the Coulomb interaction
between the electrons as a small perturbation of the hamiltonian of the
free model, which consists in a gas of non-interacting electrons
localized in the plane $(0;x,y)$ and subjected to a strong constant
magnetic field directed along z. The free hamiltonian of each electron
is, in convenient variables, the hamiltonian of an harmonic oscillator,
and the free hamiltonian of the system is the sum of all these
contributions, one for each electron. The assumption on the Coulomb
interaction allows us to deduce that, without other interactions, all
the electrons are in the lowest Landau level (ground state of the
oscillators) and the kinetic part of the complete hamiltonian (free
hamiltonian plus Coulomb interaction) can be replaced by its mean value
in the ground state of the system, so that it is replaced by a c-number.
The idea of considering only the lowest Landau level in the description
of the FQHE is now very common in the literature and was originally
proposed by Yoshioka and Fukuyama in \cite{yosh}.

Before describing our treatment of the dynamics it is worthwhile to
remark the likeness of this model with the dissipative laser model
introduced in \cite{alsew}. Of course, the present model is quite a
simplified version of that one, first for its conservative nature, which
allows an hamiltonian approach, second for the presence of only one mode
for the bosons (condition which will be relaxed in the next subsection),
and finally for the simplified form of the interaction hamiltonian. As
it has been shown in \cite{alsew, bagsew}, a complete analysis of the
Alli-Sewell model can be carried out already in a non-topological
*-algebra of unbounded operators. In this paper we want to play a similar
game with rather different rules, that is to introduce a topology somehow
related to the model itself.

The way of
doing this is suggested from the results of Section 3 and by the results
in \cite{bm}: the algebra $\A$ should reflect the composite nature of the
model, and is therefore reasonable to take $\A$ as the tensor product
of the bounded spin operators, $B(\Hil_{spin})$, and of the Lassner
boson algebra, $\Lc^+(\D)$: $\A= B(\Hil_{spin}) \otimes \Lc^+(\D)$. The
topology on $\A$, $\tau$, is generated by seminorms which deal
separately with the spin and the boson variables.  Obviously, in view of
the results in Section 3 for the free bosons and what is known about the
mean field spin models, see \cite{bm}, we take these seminorms
essentially as the ones in (\ref{27}) for the boson operators and
as the strong ones in (\ref{29}) for the spin observables. We give here
the definition of  these seminorms  which will be discussed in some more
details in Appendix B. First of all, following \cite{bm,bt}, we give the
definition of what we will call {\em relevant } vectors in the present
contest. It may be worthwhile to observe that the set below is slightly
different from the one introduced in Section 2:   
\be  \F =\left\{ \Psi \in \Hil_{spin} \: :
\lim_{|V|,\infty} \frac{1}{|V|}\sum_{p\in V} {\sigma_3^p}\Psi =
\sigma_3^\infty \Psi, \:\: \|\sigma_3^\infty\|\leq1 \right\},
\label{43}
\en
$\sigma_3^\infty$  being a certain element of the center of the
$C^*$-spin algebra. The usefulness of this set relies on the fact that
$\sigma_3^V$ (together with its powers and its analytical functions) only
converges strongly on the vectors of $\F$, not on general vectors of
$\Hil_{spin}$ (and, of course, not in norm). 

For future convenience, it is better to
introduce the unbounded operator \be M\equiv N+\Id=a^\dagger a
+\Id=\sum_{l=1}^{\infty}l\Pi_l. \label{44}
\en
As we can see, the operators $M$ and $N$ are essentially the same
observable. The only difference is in the lower value of $l$ in their
spectral decompositions. With these considerations we define, for each
$X\in \Lc^+(\D)$ and for each $A\in B(\Hil_{spin})$, 
\be
\|XA\|^{f,k,\Psi}\equiv \|X\|^{f,k} \|A\|^\Psi=\|f(M)XM^k\| \|A\Psi\|,
\label{45}
\en
where $f\in \C$, $k\in \N$, $\Psi \in \F$. The final remark is that,
again, we are considering only one contribution in the Lassner seminorm;
a complete, but unnecessary here, definition should include also
$\|M^kXf(M)\|$, see (\ref{24}).

Now we proceed to a complete regularization of the hamiltonian $H_V$.
The approach we follows is the one already discussed in Section 3.
Therefore, we consider the following operator
\be
H_{V,L}=J|V|(\sigma_3^V)^2 + \gamma (a_L+a_L^\dagger)\sigma_3^V,
\label{46}
\en
where, as before, $a_L=Q_LaQ_L=aQ_L$, see (\ref{a2}). Then we define
the cutoffed algebraic dynamics $\alpha_{V,L}^t$ in the canonical way:
\be
\alpha_{V,L}^t(X)=e^{iH_{V,L}t}Xe^{-iH_{V,L}t},
\label{47}
\en
where $X\in \A$.

For technical convenience it may be useful to introduce here a relation
between the two different cutoffs. In the second part of Proposition 5.1
we will assume the following relation between $|V|$ and $L$:
\be
|V|=L^r,
\label{48}
\en
where $r$ is an integer bigger than $1$. This condition
essentially fix the subspaces $V\subset \Z^d$ which are to be used in the
regularization procedure: they are those finite lattices whose number of
sites is the $r$-th power of the integers.

We can now state the main result of this Section:

\vspace{4mm}

{\bf Proposition 5.1}-- The limit of $\alpha_{V,L}^t(a)$ for $|V|$ and
$L$ both diverging exists in $\A[\tau]$.
In the hypothesis (\ref{48}) the same holds true also for
$\alpha_{V,L}^t(\sigma_\alpha^i)$.

\vspace{3mm}

{\bf Proof}

We begin with the proof of the first statement above. Since
$[J|V|(\sigma_3^V)^2, \gamma (a_L+a_L^\dagger)\sigma_3^V]]=0$ we can write
\bea
\alpha_{V,L}^t(a)&=&e^{i\gamma
(a_L+a_L^\dagger)\sigma_3^Vt}\left(e^{iJ|V|(\sigma_3^V)^2t}a
e^{-iJ|V|(\sigma_3^V)^2t}\right)e^{-i\gamma
(a_L+a_L^\dagger)\sigma_3^Vt}=\nonumber\\
&=&e^{i\gamma (a_L+a_L^\dagger)\sigma_3^Vt}a e^{-i\gamma
(a_L+a_L^\dagger)\sigma_3^Vt}=\sum_{n=0}^\infty\frac{(i\gamma
t\sigma_3^V)^n}{n!}[a_L+a_L^\dagger,a]_n, 
\label{49}
\ena
due to the commutation relation $[a,\sigma_\alpha^i]=0$. Here
$[a_L+a_L^\dagger,a]_n$ is the usual multiple commutator defined as
$$
[a_L+a_L^\dagger,a]_0=a, \quad \quad
[a_L+a_L^\dagger,a]_n=[(a_L+a_L^\dagger),[a_L+a_L^\dagger,a]_{n-1}], \:
n\geq 1.
$$
We divide the infinite sum in the rhs of (\ref{49}) as 
\be
\alpha_{V,L}^t(a)=\sum_{n=0}^1\frac{(i\gamma
t\sigma_3^V)^n}{n!}[a_L+a_L^\dagger,a]_n + \sum_{n=2}^\infty\frac{(i\gamma
t\sigma_3^V)^n}{n!}[a_L+a_L^\dagger,a]_n,
\label{410}
\en
and we show that the last contribution, $\sum_{n=2}^\infty...$, converges
to zero in $\tau$ when $L$ go to infinity, for any $|V|$.

The proof of this convergence goes as follows:

(a) first we compute, using the formulas in Appendix A, the following
multiple commutators:
\bea
&&[a_L+a_L^\dagger,a]_0 = a, \hspace{1cm}
[a_L+a_L^\dagger,a]_1=a^2\Pi_{L+1}+aa^\dagger\Pi_L-Q_L,\nonumber\\
&&[a_L+a_L^\dagger,a]_2 = a^3\Pi_{L+1}+a^2a^\dagger(\Pi_L+\Pi_{L+1})
-2a\Pi_{L+1}-a(a^\dagger)^2\Pi_{L-1}. 
\label{411}
\ena

(b) then we sketch the proof of the following fact:
$\|f(M)[a_L+a_L^\dagger,a]_2M^k\|\rightarrow 0$ when $L\rightarrow
\infty$, for any $f\in \C$ and any $k\in \N$. The proof is a bit long.
Here we only show the convergence to zero of the first contribution in
(\ref{411}). We have  $$ \|f(M)a^3\Pi_{L+1}M^k\|=\sum_{l,s=1}^\infty
f(l)s^k\|\Pi_la^3\Pi_{L+1}\Pi_{s}\|=
\sum_{l,s=1}^\infty f(l)s^k\delta_{s,L+1}\|\Pi_la^3\Pi_{s}\|. 
$$
But, since $\|\Pi_la\Pi_s\|=\sqrt{s}\delta_{l,s-1}$, we deduce that, using
the idempotence of $\Pi_l$, $\Pi_l=\Pi_l^2$ together with formula
(\ref{a2}),  
$$
\|\Pi_la^3\Pi_{s}\|\leq \|\Pi_la^3\Pi_{l+1}\|\|\Pi_{l+1}a^3\Pi_{l+2}\|
\|\Pi_{l+2}a^3\Pi_{s}\|\leq \sqrt{(l+1)(l+2)(l+3)}\delta_{s,l+3}.
$$
Therefore, after some algebraic computations,
$$
\|f(M)a^3\Pi_{L+1}M^k\|\leq f(L-2)(L+1)^k\sqrt{L(L-1)(L-2)},
$$
which, of course, goes to zero for $L\rightarrow \infty$, due to decaying
features of the function $f$. Analogous estimates work for the other
contributions in $\|f(M)[a_L+a_L^\dagger,a]_2M^k\|$ in (\ref{411}).

(c) in the third step we show how the estimate before can be used also
for the other multiple commutators $[a_L+a_L^\dagger,a]_n$, $n>2$. The estimate goes
in the following way: first we observe that $[a_L+a_L^\dagger,a]_2$ is
essentially the sum of 6 monomials (We consider two contributions from
$-2a\Pi_{L+1}$) in $a$ and $a^\dagger$, each at most raised to
the power 3. Since we are interest to the limit of very large $L$, it is
not very important the fact that sometimes we find $\Pi_L$ and somewhere
else $\Pi_{L\pm1}$. The structure of this commutator is therefore
$[a_L+a_L^\dagger,a]_2 \simeq 6 m_2(a,a^\dagger)\Pi_L$, where $m_2$ is a
monomial of degree $2+1=3$.
Analogously we find $[a_L+a_L^\dagger,a]_3 \simeq 24
m_3(a,a^\dagger)\Pi_L$ where $m_3$ is of degree 4. The number 24 is the
maximum possible number of contributions following from the commutation
between $a_L+a_L^\dagger$ and $6 m_2(a,a^\dagger)\Pi_L$. This procedure
can be extended  to each $n$, and we find that, for $L \gg 1$,
\be
[a_L+a_L^\dagger,a]_n \simeq 6\cdot 4^{n-2} m_n(a,a^\dagger)\Pi_L.
\label{412}
\en
Using this formula in the expression of
$\|[a_L+a_L^\dagger,a]_n\|^{f,k}$, together with the
$\Pi_L\Pi_s=\Pi_s\delta_{s,L}$ we obtain the following:
$$
\|[a_L+a_L^\dagger,a]_n\|^{f,k}\simeq 6L^k4^{n-2}\sum_{l=1}^\infty
\|\Pi_lm_n(a,a^\dagger)\Pi_L\|.
$$
Using now the same kind of estimates used in point (b) above, we see
that  $\|\Pi_lm_n(a,a^\dagger)\Pi_L\|$ is essentially a
polinomial $p_n$ in $l$ and $L$, times a delta function
$\delta_{l,\varphi_{m_n}}$, where $\varphi_{m_n}$ is an increasing
function related to the explicit expression of $m_n$. Therefore we get
$$
 \|[a_L+a_L^\dagger,a]_n\|^{f,k}\simeq 6 \cdot 4^{n-2} L^k
f(\varphi_{m_n}(L))p_n(\varphi_{m_n}(L),L).
$$
Due to the nature of $f$, the right hand side is bounded by the supremum
on $n$ of the same quantity:
$$
L^k f(\varphi_{m_n}(L))p_n(\varphi_{m_n}(L),L) \leq \sup_n 
{\Large[}L^k f(\varphi_{m_n}(L))p_n(\varphi_{m_n}(L),L){\Large]} =:F_k(L),
$$
where the function $F_k(L)$ surely goes to zero when $L\rightarrow
\infty$. We conclude that
\be
\|[a_L+a_L^\dagger,a]_n\|^{f,k} \leq  6 \cdot 4^{n-2} F_k(L).
\label{413}
\en

\vspace{4mm}

(d) At this point we have at hand all the ingredients to prove that the
last contribution in (\ref{410}) converges to zero for any $V$ when $L$
diverges. In fact we have, using the definition (\ref{45}) of the
seminorms $\|\:\|^{f,k,\Psi}$, the inequality (\ref{413}), and the
uniform bound on $\sigma_3^V$, $\|(\sigma_3^V)^n \| \leq 1$, 
\beano
& &\|\sum_{n=2}^\infty\frac{(i\gamma
t\sigma_3^V)^n}{n!}[a_L+a_L^\dagger,a]_n \|^{f,k,\Psi} \leq 
\sum_{n=2}^\infty \frac{|\gamma t|^n}{n!} \|(\sigma_3^V)^n
\Psi\|\|[a_L+a_L^\dagger,a]_n \|^{f,k}\leq \\
& & \frac{|\gamma t|^2}{2!}\|[a_L+a_L^\dagger,a]_2 \|^{f,k}+ 
\frac{|\gamma t|^3}{3!}\|[a_L+a_L^\dagger,a]_2 \|^{f_1,k_1}+
\frac{|\gamma t|^4}{4!}\|[a_L+a_L^\dagger,a]_2 \|^{f_2,k_2}+...\leq \\
& & \leq \frac{3}{2}\left(\sum_{n=2}^\infty \frac{|4\gamma t
|^n}{n!}\right) F_k(L) \leq \frac{3}{2}e^{|4t\gamma|}F_k(L)
\rightarrow 0,
\enano
when $L\rightarrow \infty$.

Now the first statement of the
Proposition, that is the $\tau$-convergence of $\alpha_{V,L}^t(a)$,
easily follows.
 In fact, in view of the above result, we can write equality (\ref{410})
as \be
\alpha_{V,L}^t(a)=a+i\gamma \sigma_3^Vt(a^2\Pi_{L+1}+aa^\dagger\Pi_L-Q_L)
+\tilde\alpha_{V,L}^t(a),
\label{417}
\en
where $\tilde\alpha_{V,L}^t(a)\equiv \sum_{n=2}^\infty\frac{(i\gamma
t\sigma_3^V)^n}{n!}[a_L+a_L^\dagger,a]_n$ is $\tau$-converging to zero with
$L$, and uniformly bounded in $|V|$, so that 
\beano
& &\||\alpha_{V,L}^t(a)-\alpha_{V',L'}^t(a)\||^{f,k,\Psi} \leq \\
& &\leq |\gamma t| 
\||\sigma_3^Vt(a^2\Pi_{L+1}+aa^\dagger\Pi_L-Q_L)-\sigma_3^{V'}t(a^2\Pi_{L'+1}+
aa^\dagger\Pi_{L'}-Q_{L'})\||^{f,k,\Psi}+ O_L,
\enano
where $O_L\rightarrow 0$ when $L\rightarrow \infty$.

In this estimate we can handle very easily all the contributions
containing any projection operator $\Pi_n$, such as
$\||\sigma_3^Vta^2\Pi_{L+1}\||^{f,k,\Psi}$. In fact, they go to zero, how
can be seen following essentially the same procedure as in
point (b) above. What is to be controlled is only the term
$\||\sigma_3^VQ_L- \sigma_3^{V'}Q_{L'}\||^{f,k,\Psi}$. This can be
estimated adding and subtracting first the same quantity
$\sigma_3^VQ_{L'}$ and observing that, if $L>L'$, then $$
\|Q_L-Q_{L'}\|^{f,k}=\sum_{l=L'+1}^Lf(l)l^k\rightarrow 0, $$ when $L$ and
$L'$ both diverge. We have  \beano &
&\||\sigma_3^VQ_L-\sigma_3^{V'}Q_{L'}\||^{f,k,\Psi} \leq 
\||\sigma_3^V(Q_L-Q_{L'})\||^{f,k,\Psi} +
\||(\sigma_3^V-\sigma_3^{V'})Q_{L'}\||^{f,k,\Psi} = \\ & & =
\|\sigma_3^V\Psi\|\:
\|Q_L-Q_{L'}\|^{f,k}+\|(\sigma_3^V-\sigma_3^{V'})\Psi\|\:
\|Q_{L'}\|^{f,k}\rightarrow O \hspace{4mm} \mbox{for } L,|V|\rightarrow
\infty. \enano The existence of $\alpha^t(a)$ finally follows from the
completeness of the algebra $\A$.

\vspace{4mm}

{\bf Remarks.}-- (1) It is clear from the above procedure that the order
of the limits ($|V|\rightarrow \infty$ or $L\rightarrow \infty$) has no
importance.

(2) It is worthwhile to observe that we have not used yet the
hypothesis (\ref{48}). We will use it in the second part of the proof.

\vspace{6mm}

Let us now start with the proof of the second statement of the
Proposition, that is the convergence of $\alpha_{V,L}^t(\sigma_\alpha^i)$.
Using the commutativity of the two pieces of the hamiltonian,
formula (2.3) of reference \cite{bt}, and the commutation relation
$[\sigma_3^V,\sigma_3^i]=0$, we have 
\bea
&&\hspace{-17mm}\alpha_{V,L}^t(\sigma_\alpha^i)=e^{iH_{V,L}t}\sigma_\alpha^i
e^{-iH_{V,L}t}= \nonumber\\ 
&&\hspace{-17mm}=e^{i\gamma
(a_L+a_L^\dagger)\sigma_3^Vt}\left(e^{iJ|V|(\sigma_3^V)^2t}
\sigma_\alpha^i e^{-iJ|V|(\sigma_3^V)^2t}\right)e^{-i\gamma
(a_L+a_L^\dagger)\sigma_3^Vt}=\beta_{V,L}^t(\sigma_\alpha^i)\cos^2(S_3^V)-
\\
&&\hspace{-17mm}-2\epsilon_{3\alpha\beta}
\beta_{V,L}^t(\sigma_\beta^i)\sin(S_3^V)\cos(S_3^V)+\sigma_3^i \,
\beta_{V,L}^t(\sigma_\alpha^i)\,\sigma_3^i\sin^2(S_3^V)+O(|V|^{-1}),
\nonumber
\label{418} 
\ena 
where $S_3^V=2Jt\sigma_3^V$, $O(|V|^{-1})$ is norm converging to zero and
we have introduced the notation 
\be
\beta_{V,L}^t(\sigma_\alpha^i)=e^{i\gamma (a_L+a_L^\dagger)\sigma_3^Vt}
\sigma_\alpha^i e^{-i\gamma (a_L+a_L^\dagger)\sigma_3^Vt}.
\label{419}
\en
The next step consists in proving that 
\be
\beta_{V,L}^t(\sigma_\alpha^i)=\sigma_\alpha^i+
\tilde\beta_{V,L}^t, 
\label{420}
\en
where $\tilde\beta_{V,L}^t$ 
$\tau$-converges to zero when $|V|$ diverges under the technical
assumption (\ref{48}). This result is reasonable, since
$[\sigma_3^V,\sigma_\alpha^i]=O(|V|^{-1})$, so that in definition
(\ref{419}) we can commute $e^{i\gamma (a_L+a_L^\dagger)\sigma_3^Vt}$
with $\sigma_\alpha^i$ only within an error going to zero when $|V|$ is
sent to infinity. We are going to prove (\ref{420}) rigorously. We have 
\bea
 \beta_{V,L}^t(\sigma_\alpha^i)&=&\sum_{l=0}^\infty
\frac{(it\gamma(a_L+a_L^\dagger))^l}{l!}[\sigma_3^V,\sigma_\alpha^i]_l=
\nonumber\\ &=& \sigma_\alpha^i+ \sum_{l=1}^\infty
\frac{(it\gamma(a_L+a_L^\dagger))^l}{l!}[\sigma_3^V,\sigma_\alpha^i]_l.
\label{421} \ena
We will show that the last term, which is exactly what we have called
$\tilde\beta_{V,L}^t$, converges to zero in $\tau$, at least under the
assumption (\ref{48}). Since $$
\||\sum_{l=1}^\infty \frac{(it\gamma(a_L+a_L^\dagger))^l}{l!}
[\sigma_3^V,\sigma_\alpha^i]_l\||^{f,k,\Psi} \leq \sum_{l=1}^\infty
\frac{|t\gamma|^l}{l!}\|(a_L+a_L^\dagger)^l\|^{f,k} 
\|[\sigma_3^V,\sigma_\alpha^i]_l\Psi\|,
$$
what we need to estimate are $\|(a_L+a_L^\dagger)^l\|^{f,k}$ and 
$\|[\sigma_3^V,\sigma_\alpha^i]_l\Psi\|$. This last term is easily seen
to be bounded by powers of $2/|V|$:
\be
\|[\sigma_3^V,\sigma_\alpha^i]_l\Psi\| \leq \left(\frac{2}{|V|}\right)^l.
\label{422}
\en
For what concerns the estimate of the $\|(a_L+a_L^\dagger)^l\|^{f,k}$ the
situation is a bit more involved. In fact, it is possible to prove the
following inequality
\be
\|(a_L+a_L^\dagger)^l\|^{f,k} \leq c_k \left((2^k+1)(L+1)^{3/2}\right)^l
\label{423}
\en
where $c_k\equiv \sum_{l=1}^\infty f(l)l^k$ is certainly finite for any
integer $k$. The proof of the estimate above goes as follows:
$$
\|(a_L+a_L^\dagger)^l\|^{f,k}=\|f(M)(a_L+a_L^\dagger)^lM^k\|\leq
\|f(M)M^k\| \|M^{-k}(a_L+a_L^\dagger)M^k\|^l.
$$
Now, $\|f(M)M^k\|=c_k$ and, using the spectral decomposition for $M$ and
formula (\ref{a4}), 
$$
\|M^{-k}(a_L+a_L^\dagger)M^k\|= \sum_{l,s=1}^L\left(\frac{s}{l}\right)^k
\|\Pi_l(a+a^\dagger)\Pi_s\| \leq (2^k+1)L\sqrt{L+1}.
$$
Putting all together inequality (\ref{423}) follows. Now we can use this
inequality and the one in (\ref{422}) to prove our claim on 
$\tilde\beta_{V,L}^t$. In fact we have
\beano
&&\||\sum_{l=1}^\infty \frac{(it\gamma(a_L+a_L^\dagger))^l}{l!}
[\sigma_3^V,\sigma_\alpha^i]_l\||^{f,k,\Psi} \leq c_k \sum_{l=1}^\infty
\left(\frac{(|2t\gamma|(2^k+1)(L+1)^{3/2})}{|V|}\right)^l\frac{1}{l!}=\\
&& = c_k\left(
\exp\left\{\frac{(|2t\gamma|(2^k+1)(L+1)^{3/2})}{|V|}\right\}-1\right),
\enano
which, of course, goes to zero for $L\rightarrow \infty$, at least if 
the hypothesis (\ref{48}) is satisfied. This concludes the proof of the
equation (\ref{420}). Taking into account condition (\ref{48}) we will
write  \be
\||\beta_{L}^t(\sigma_\alpha^i)-\sigma_\alpha^i\||^{f,k,\Psi}
=O(L^{-\mu}), 
\label{424}
\en
where $\mu$ is a strictly positive number. Inserting (\ref{424}) in
(\ref{418}), with $|V|=L^r$, we get
\be
\alpha_{L}^t(\sigma_\alpha^i)=\sigma_\alpha^i\cos^2(S_3^{L^r})-
2\epsilon_{3\alpha\beta} \sigma_\beta^i\sin(S_3^{L^r})\cos(S_3^{L^r})+
\sigma_3^i
\sigma_\alpha^i\sigma_3^i\sin^2(S_3^{L^r})+O(L^{-r})+O(L^{-\mu}).
\label{425}
\en
Here we have written explicitly the two terms going to zero
with $L$, the reason being that while the first  goes to zero in the usual
norm of the operators, see \cite{bt}, the second one goes to zero in the
topology $\tau$. Using now the expression above for
$\alpha_{L}^t(\sigma_\alpha^i)$ and the results in \cite{bm}, we conclude
that the sequence $\alpha_{L}^t(\sigma_\alpha^i)$ is $\tau$-Cauchy, so
that its limit exists in the algebra $\A$.
\hfill$\Box$

\vspace{4mm}

{\bf Remarks:}-- (a) It is worthwhile to notice that condition (\ref{48}) is
only a technical tool to simplify the computations. We do believe that it
is not an essential assumption.

(b) It is also worthwhile to remark that, even if the approach we
are following in this paper is somehow more in the line of \cite{bt,
bt2}, nevertheless is reference \cite{bm} which gives us the
information we need about what we have called '$\F$-strong convergence'
of $\sigma_3^V$ and of its {\em good}, i.e. analytical, functions. In
particular, in \cite{bt} we are forced not to use the strong topology
because we deal with variables, $\frac{1}{|V|^\gamma}\sum_{i\in
V}\sigma_3^i$, which are, in general, not uniformly bounded in $V$, since
we only assume $\gamma \leq 1$. Of course, the algebraic structure for
the spin variables of this model, could also be constructed using the
results in \cite{bt} which, for $\gamma=1$, is nothing but a mean field
model, like the present model. This would imply, however, that the strong
seminorms in (\ref{29}) should be replaced with the Lassner seminorms
constructed starting with a certain number operator which can be
introduced also in the spin algebra, \cite{las,bt}. Also, the algebra
$B(\Hil_{spin})$ should be replaced by a certain $\Lc^+(\tilde\D)$, where
$\tilde\D$ is the dominion of all the powers of this new number operator,
and is a subspace of $\Hil_{spin}$. Since all of this is not necessary
here, being $\sigma_3^V$ uniformly bounded, we prefer to use the simpler
framework discussed in \cite{bm}.

(c) With the same kind of estimates we have used above, it is also
possible to prove the existence of, among the others, the 
limits of the following quantities: $\alpha_L^t(a\sigma_\alpha^i)$,
$\alpha_L^t(\sigma_{\alpha_1}^{i_1}.....\sigma_{\alpha_n}^{i_n})$,
$\alpha_L^t(a^n)$. Here $n$ is an arbitrary integer. The proof
essentially follows from the explicit expressions (\ref{417}) and
(\ref{425}) for the cutoffed dynamics and from the algebraic nature of
$\alpha_L^t$.

\vspace{4mm}

We conclude, therefore, that the algebraic dynamics for the model in
(\ref{42}) can be defined using the regularization (\ref{34bis}) and
that a natural topological *-algebra where to consider the model is
$\A[\tau]=B(\Hil_{spin})\otimes \Lc^+(\D)$, $B(\Hil_{spin})$ taken with
its strong topology and $\Lc^+(\D)$ with the physical topology.

%4.2
\subsection{\hspace{-6mm}. Many mode bosons}

The last model we will discuss in this paper is an extension of the
interacting model described by the hamiltonian in (\ref{41}):
\be
H'_V=\frac{J}{|V|}\sum_{i,j\in
V}\sigma_3^i\sigma_3^j+\sum_{l=1}^n a^\dagger_l a_l +\sum_{l=1}^n\gamma_l
(a_l+a_l^\dagger)\sigma_3^V.
\label{426}
\en
In this generalization, $n$ different modes of bosons are introduced,
satisfying the CCR $$
[a_l,a^\dagger_n]=\delta_{ln}\Id,
$$
while $\gamma_l$ are $n$ small real numbers such that, again, the
kinetic boson hamiltonian $\sum_{l=1}^n a^\dagger_l a_l$ can be
considered frozen in a given level. Again, therefore, we consider the
following hamiltonian \be
H_V=\frac{J}{|V|}\sum_{i,j\in
V}\sigma_3^i\sigma_3^j +\sum_{l=1}^n\gamma_l
(a_l+a_l^\dagger)\sigma_3^V.
\label{427}
\en
The regularization of $H_V$ can be performed in the following way: for
each mode we consider the number operator $N_k=a_k^\dagger a_k$ and its
spectral decomposition $N_k=\sum_{l=0}^\infty l \Pi_{k,l}$ and the
projection operators $Q_{k,L}\equiv \sum_{l=0}^L \Pi_{k,l}$ (Also for
this model it is convenient to introduce $M_k=N_k+\Id$). We also
define the Hilbert spaces $\Hil_k$, one for each mode, and $\D_k \equiv
D^\infty(N_k) \subset \Hil_k$ as in Section 2. Finally we take $\D\equiv
\D_1\otimes\D_2\otimes.....\otimes\D_n$ which is dense in the Hilbert
space $\Hil$ of the tensor product of all the $\Hil_k$. More useful
spaces are $\E_{k,l}=\{(a_k^\dagger)^l\Phi_o\}$ and $\F_{k,L}=
\E_{k,0}+\E_{k,1}+......+\E_{k,L}$. Of course all the above projection
operators related to different modes commute. The algebra is the usual
$\A=B(\Hil_{spin})\otimes\Lc^+(\D)$.

At this point we define the regularized hamiltonian $H_{V,L}$ as
\be
H_{V,L}=\frac{J}{|V|}\sum_{i,j\in
V}\sigma_3^i\sigma_3^j +\sum_{l=1}^n\gamma_{l}
(a_{l,L}+a_{l,L}^\dagger)\sigma_3^V,
\label{428}
\en
where $a_{l,L}=Q_{l,L}a_lQ_{l,L}$.

The existence of the limit of $\alpha_{V,L}^t(a_j)$ is a simple
consequence of the results in the first subsection. In fact, with the
same arguments used for the single mode bosons model, we find that 
$$
\alpha_{V,L}^t(a_j)=e^{i\gamma (a_{j,L}+a_{j,L}^\dagger)\sigma_3^Vt}a_j
e^{-i\gamma (a_{j,L}+a_{j,L}^\dagger)\sigma_3^Vt}.
$$
This has exactly the same form seen in equation (\ref{49}), so that its
convergence in the physical topology can be proven in the same way.

A bit different is the proof of the existence of the limit of
$\alpha_{V,L}^t(\sigma_\alpha^i)$. Without going too much into details,
we just observe that the main difference is that in formula (\ref{418})
$\beta_{V,L}^t(\sigma_\alpha^i)$ must be replaced with
\beano
&&\varphi_{V,L}^t(\sigma_\alpha^i)=e^{i\sum_{l=1}^n\gamma_l
(a_{l,L}+a_{l,L}^\dagger)\sigma_3^Vt} \sigma_\alpha^i e^{-i
\sum_{l=1}^n\gamma_l (a_{l,L}+a_{l,L}^\dagger)\sigma_3^Vt}=\\
&=&e^{i\gamma_n(a_{n,L}+a_{n,L}^\dagger)\sigma_3^Vt}\left(....
\left(e^{i\gamma_1(a_{1,L}+a_{1,L}^\dagger)\sigma_3^Vt} \sigma_\alpha^i
e^{-i\gamma_1(a_{1,L}+a_{1,L}^\dagger)\sigma_3^Vt}\right)....\right)
e^{-i\gamma_n(a_{n,L}+a_{n,L}^\dagger)\sigma_3^Vt}.
\enano
This can still be estimated as in (\ref{420}), observing that the
contribution $\tau$ converging to zero appears now $n$
times, but, since $n$ is the finite number of modes of our model,
formula (\ref{425}) still holds in this case, so that its limit exists
in $\A$.

% Section 5
\section{Outcome and Future Projects}

In this paper some physically relevant models have been analyzed
making use of the algebras of unbounded operators of the Lassner type.
In particular, we have shown that for the free bosons, and for other
models related to this one, it is possible to introduce an "occupation
number" cutoff, whose removal can be performed in the complete
topological *-algebra $\Lc^+(\D)$, $\D$ being the domain of all the
powers of the number operator $N$.

Moreover, for an interacting model of spins and bosons, the same kind of
regularization allows a rigorous definition of the cutoffed hamiltonian.
Again, this cutoff can be removed, together with its other volume
cutoff, working in a large topological *-algebra $\A[\tau]=
B(\Hil_{spin})\otimes \Lc^+(\D)$, where $B(\Hil_{spin})$  and $\Lc^+(\D)$
are endowed respectively with the strong and with the physical
topologies.

\vspace{3mm}

For what concerns our future projects, our main goal is to enrich the
present analysis with some more physically realistic models like the
original conservative Dicke model, \cite{dic}, or its non-conservative
generalization proposed by Alli and Sewell, \cite{alsew}.

\vspace{6mm}

\noindent{\large \bf Acknowledgments} \vspace{3mm}

	It is a pleasure to thank Prof. G.L. Sewell for many stimulating
discussions about the interactions between matter and radiation and for
his comments. I also would like to thank Dr. D. Dubin for an interesting
discussion concerning topological algebras of unbounded operators, and
Dr. C. Trapani for his careful reading of the manuscript and for his
comments. This work has been supported by M.U.R.S.T.

\vspace{8mm}

\appendix
\renewcommand{\theequation}{\Alph{section}.\arabic{equation}}

% Section 1

 \section{\hspace{-5mm} Appendix:  Some Useful Formulas}

In this Appendix we will obtain some useful relations concerning the
projection operators $\Pi_l$ and $Q_L$. We start by proving the following
relation:
\be
\Pi_{l-1}a=a\Pi_l, \hspace{1cm} l=1,2,3,...
\label{a1}
\en
The proof is the following. First, we observe that $\varphi\in \Hil$ is
such that $ \Pi_{l-1}a\varphi=0$ if, and only if, $a\Pi_{l}\varphi=0$. In
fact, if $\Pi_{l-1}a\varphi=0$ then $a\varphi$ must be orthogonal to
$\E_{l-1}$ or, that is the same, $\varphi$ must be orthogonal to $\E_l$,
which implies that $\Pi_l\varphi=0$. The converse implication is proved
in the same way. This implies that if $\Psi$ is such that $ \Pi_{l-1}a\Psi
\neq0$ then necessarily also $a\Pi_{l}\Psi \neq 0$. Let $\Psi$ be such a
vector. Then $\Pi_l \Psi \in \E_l$ and $a\Pi_l \Psi \in \E_{l-1}$.
Analogously $\Pi_{l-1}a\Psi \in \E_{l-1}$. Therefore, the two vectors must
be necessarily proportional to each other. It must exists, in other words,
a non zero constant $\alpha_l$ such that $\Pi_{l-1}a\Psi = \alpha_l a\Pi_l
\Psi$. Taking the scalar product of both sides of this equality with the
vector $(a^\dagger)^{l-1}\Phi_0$ we deduce that $\alpha_l =1$, from which
equation (\ref{a1}) follows. 

The following relations are now easy consequences of formula (\ref{a1}):
\be
a^\dagger \Pi_{l-1}=\Pi_l a^\dagger
\hspace{1cm}Q_La=aQ_{L+1}, \hspace{.5cm} Q_La^\dagger = a^\dagger Q_{L-1},
\label{a2}
\en
\be
[Q_L,a]=\Pi_La, \hspace{1cm} [Q_L,a^\dagger]=-a^\dagger \Pi_L,
\label{a3}
\en
where $l$ and $L$ take integer values.

Let us finally prove the following useful formula: 
\be
\|\Pi_la\Pi_s\|=\sqrt{s}\delta_{l,s-1}.
\label{a4}
\en
This is a consequence of the following equalities:
\beano
\|\Pi_la\Pi_s\|^2&=&\|a\Pi_{l+1}\Pi_s\|^2=\delta_{l,s-1}\|a\Pi_s\|^2= 
\delta_{l,s-1}\sup_{\|\varphi\|\leq1}|<a\Pi_s\varphi,a\Pi_s\varphi>|=\\
&=&\delta_{l,s-1}\sup_{\|\varphi\|\leq1}|<\Pi_s\varphi,N\Pi_s\varphi>|= 
s\delta_{l,s-1}.
\enano

\vspace{4mm}

%\appendix
%
%\renewcommand{\theequation}{\Alph{section}.\arabic{equation}}

% \Section 2

 \section{\hspace{-5mm} Appendix:  The topology $\tau$ for the
Interacting Models}

Let $\Hil$ and $\Hil_{spin}$ be the Hilbert spaces respectively of the
bosons and the spin. $\D$ is the subset of $\Hil$ defined as in
Section 2 and $\A=B(\Hil_{spin})\otimes \Lc^+(\D)$ is the relevant
*-algebra for the model in Section 4. 

We start defining, for $X=X_1\otimes X_2\in B(\Hil_{spin})\otimes
B(\Hil)$, and for $\Psi \in \F$ given in (\ref{43}), the following
'strong' seminorms: \be
\|X\|^\Psi=\sup_{\|\chi_1\|\leq 1,\, \|\chi_2\|\leq 1,\, \|\phi_2\|\leq 1}
|<\chi_1\otimes\chi_2, (X_1\otimes X_2) \Psi\otimes\phi_2>|= \|X_1\Psi\|
\|X_2\|,
\label{b1}
\en
(obviously in $\|X_1\Psi\|$ the norm is the one in $\Hil_{spin}$, while
in $\|X_2\|$ is the one in $\Hil$). It is straightforwardly proven that
these are really seminorms. It is also evident the meaning of the above
definition: we are considering the $(\F-)$strong topology on
$B(\Hil_{spin})$ and the usual norm topology on $B(\Hil)$. At this point,
since if $A\in \Lc^+(\D)$ then both $f(M)AM^k$ and $M^kAf(M)$ are bounded
operators on $\Hil$, the following definition appears as the most
natural: for any $X\in \A$ we define the following seminorms
\be
\||X\||^{f,k,\Psi}\equiv \max \left\{\|f(M)XM^k\|^\Psi,\|M^kXf(M)\|^\Psi
\right\}.
\label{b2}
\en
Again, it is not difficult to prove that $\{\||.... \||^{f,k,\Psi}\}$
is a system of seminorms which define a topology $\tau$ and this
topology makes $\A$ a complete topological *-algebras.

\end{document}